\documentclass[epic,eepic,12pt]{amsart}

\usepackage{amsfonts,mathrsfs}
\usepackage[mathscr]{eucal}

\usepackage{amssymb}
\usepackage{amscd}
\usepackage{fancyhdr}
\usepackage{euler,eucal}

\DeclareMathAlphabet{\mathbf}{T1}{ppl}{bx}{n}
\DeclareMathAlphabet{\mathrm}{T1}{ppl}{m}{n}

\numberwithin{equation}{section}

\newcommand\note[1]%
{$^\dagger$\marginpar{\footnotesize{$^\dagger${#1}}}}

\def\({\left(}
\def\){\right)}
\def\<{\left<}
\def\>{\right>}

\newtheorem{theorem}{Theorem}[section]
\newtheorem{Proposition}[theorem]{Proposition}
\newtheorem{lemma}[theorem]{Lemma}

\theoremstyle{definition}
\newtheorem{example}[theorem]{Example}
\newtheorem{remark}[theorem]{Remark}

\newcommand\bb[1]{{\text{\bf#1}}}

\newcommand\Z{\bb{Z}}
\newcommand\Q{\bb{Q}}
\newcommand\R{\mathbb{R}}

\newcommand\funclim[1]{\operatorname*{\mathrm{#1}}}

\renewcommand\lim{\funclim{lim}}

\newcommand\sur{\mathrel{\to\kern-1.8ex\to}}
\newcommand\iso{\mathrel{\hookrightarrow\kern-1.8ex\to}}

\newcommand\longhookrightarrow{\lhook\joinrel\longrightarrow}

\newcommand\longsur{\mathrel{\longrightarrow\kern-1.8ex\to}}
\newcommand\longiso{\mathrel{\longhookrightarrow\kern-1.8ex\to}}

\begin{document}

\bibliographystyle{amsalpha}
\date{\today}

\title{The log-concavity conjecture for the Duistermaat-Heckman
measure revisited}

\author{Yi Lin}
\maketitle
\begin{abstract}

 Karshon constructed the first counterexample to the log-concavity
 conjecture for the Duistermaat-Heckman measure: a Hamiltonian six
 manifold whose fixed points set is the disjoint union of two copies of $T^4$.
 In this article, for any closed symplectic four manifold $N$ with $b^+>1$, we show that there is
a Hamiltonian six manifold $M$ such that its fixed points set is the
disjoint union of two copies of $N$ and such that its
Duistermaat-Heckman function is not log-concave.

On the other hand, we  prove that if there is a torus action of
complexity two such that all the symplectic reduced spaces taken at
regular values satisfy the condition $b^+=1$, then its
Duistermaat-Heckman function has to be log-concave. As a
consequence, we prove the log-concavity conjecture for  Hamiltonian
circle actions on six manifolds such that the fixed points sets have
no four dimensional components, or only have four dimensional pieces
with $b^+=1$.

\end{abstract}

\section{Introduction}

Consider the effective Hamiltonian action of a torus $T$ on a
$2n$-dimensional  connected symplectic manifold $(M,\, \sigma)$ with
a proper moment map $\Phi: M \rightarrow \frak{t}^*$, where
$\frak{t}=\text{Lie}(T)$. The {\bf Duistermaat-Heckman measure}
\cite{DH82} on $\frak{t}^*$ is the push-forward of the Liouville
measure $\vert \beta \vert$, the one defined by the symplectic
volume form $\dfrac{1}{n!}\omega^n$, via the momentum map $\Phi$.
The Duistermaat\--Heckman measure is absolutely continuous with
respect to the Lebesgue measure, and its density function, which is well defined once the
normalization of the Lebesgue measure is declared, is said
to be the {\bf Duistermaat\--Heckman function}.

More generally, consider the Hamiltonian action of a compact
connected Lie group $G$ on the symplectic manifold $(M,\omega)$ with
a proper moment map $\Phi: M\rightarrow \frak{g}^*$, where
$\frak{g}=\text{Lie}(G)$. Let $T$ be the maximal torus of $G$ with
Lie algebra $\frak{t}$ and $W$ the Weyl group of $G$. Choose a
$W$-invariant inner product on $\frak{t}$ so as to identify it with
$\frak{t}^*$. Then we define measure $\nu$ on the positive Weyl chamber $\frak{h}_+
\subset \frak{t}^*$ by letting $\frac{1}{\vert W\vert} \nu$ be the pushforward of the measure $ \vert \beta \vert$ via the composition
$M\xrightarrow{\Phi} \frak{g}^* \rightarrow \frak{h}_+ =\frak{g}^*/G$.

 %the push-forward measure $\Phi_*\vert \beta
% \vert$ is well defined on the space of compactly supported
% $G$-invariant functions on $\frak{g}^*$, and is uniquely determined
% by its restriction, $\nu$, to  as follows: if $\rho$ is a $G$-invariant
% compactly supported function on $\frak{g}^*$ and $\overline{\rho}$
% is its restriction to $\frak{h}_+$, then $\frac{1}{\vert
% W\vert}\int_{\frak{h}_+}\overline{\rho}\,\, \nu
% =\int_{\frak{g}^*}\rho\,\, \Phi_*\vert \beta \vert.$

 A measure defined on $\R^k$ is said to be
log-concave if it is absolutely continuous with respect to the
Lebesgue measure on $\R^k$ and if the logarithm of its density function is
a concave function. The log-concavity of the Duistermaat-Heckman measure,
and, more generally, of the measure $\nu$ we described in the
previous paragraph, has been established for circle actions on four
manifolds by Y. Karshon \cite[Remark 2.19]{Ka94}, for torus actions
on compact K\"ahler manifolds by W. Graham \cite{Gr96}, and for
compact connected group actions on projective varieties, possibly
singular, by A. Okounkov \cite{okounkov1}, \cite{okounkov},
\cite{ok98}. In particular, the result of \cite{okounkov1} led V.
Ginzburg to conjecture that the Duistermaat-Heckman
measure is log-concave for any Hamiltonian torus action
on a compact symplectic manifold $M$.

 The same conjecture, independently of \cite{okounkov},
was propsed by A. Knutson. In view of the above positive results,
this conjecture seems very plausible in early nineties. However, motivated by an example of McDuff \cite{Mcduff}, Karshon
\cite{Ka96} constructed a Hamiltonian circle action on a compact six
manifold for which the Duistermaat-Heckman measure is
non-log-concave, which provides the first counterexample. For more background materials and
motivations of the log-concavity properties of the Duistermaat-Heckman measure, the interested readers
are strongly encouraged to read \cite{OK00} for an excellent
expository account.

 In a different direction, inspired by \cite{Ka96} and \cite{yan;hodge-structure-symplectic}, the author
\cite{Lin07} constructed the first counterexamples that the Hard
Lefschetz property does not survive the symplectic reduction. These
examples provide us with an infinite class of six dimensional
Hamiltonian circle manifolds which do not admit a K\"ahler
structure. Naturally, the author was led to the question of whether
the construction used in \cite{Lin07} can be adapted to produce
general examples of non-K\"ahler Hamiltonian manifolds with
non-log-concave Duistermaat-Heckman functions so as to offer a
better understanding why the log-concavity property fails in the
general symplectic category.

%Thus one might want to use them to test the log-concavity conjecture
%for the Duistermaat-Heckman measure.

%Let us just name a few such examples here: for complex K\"ahler
% surfaces, $K_3$ surfaces, and so on ( to be added later); for
% non-K\"ahler symplectic manifolds, lots of examples presented in
% \cite[Section 5,6]{Gom} and ( to be added later).

 In this article we will present a rather satisfactory answer to the above question. First,
 we prove that for any closed symplectic four manifold $N$ with $b^+>1$,
  there exists a symplectic six manifold fibred over $N$ such
that there is a Hamiltonian $S^1$ action on $M$ for which the
Duistermaat-Heckman function is non-log-concave. This provides us with a huge class of Hamiltonian manifolds with a non-log-concave
Duistermaat-Heckman function, since there are many symplectic manifolds which satisfy
$b^+>1$, c.f., \cite{Gom95} and \cite{PS00}. As an application, we construct
simply-connected six dimensional Hamiltonian circle manifolds which
satisfy the Hard Lefschetz property and which have a non-log-concave
Duistermaat-Heckman function. In particular, this shows that the
Hard Lefschetz property, unlike that of invariant K\"ahler
condition, does not imply the log-concavity conjecture.
% For instance, among K\"ahler surfaces, $K3$ surfaces and complex tori satisfy $b^+=3$. For general symplectic manifolds

Second, we give a useful cohomological condition which ensures the
Duistermaat-Heckman function of a complexity two Hamiltonian torus
action to be  log-concave. More precisely, we prove that if the
symplectic quotients taken at any regular value have $b^+=1$, then
the Duistermaat-Heckman function must be log-concave. As a result,
we establish the log-concavity conjecture for circle actions on six
manifolds such that the fixed points sets have no four dimensional
components, or have only four dimensional pieces with $b^+=1$.

Indeed, given a circle action on a six manifold which satisfies the
above assumptions, when the action is semi-free, i.e., it is free on the complement of
the fixed point set, the fact that the
symplectic quotients taken at regular values have $b^+=1$ can be
seen by the following observations. First, applying the equivariant
Darboux theorem to an invariant open neighborhood of the minimal
critical submanifold, one checks easily that $b^+=1$ for symplectic
quotients taken at a regular value sufficiently close to the
minimum. \footnote{ Details are given in the proof of Theorem \ref{log-concavity-conjecture}. } By the Duistermaat-Heckman theorem \cite{DH82}, in the same
connected component of the regular values of the moment map, the
diffeotype of symplectic quotients does not change. When passing a
critical level of the moment map, the diffeotype of symplectic
quotients changes by a blow up followed by a blow down \cite{GS89}.
Since the symplectic
quotients under consideration here are all four dimensional, blowing up along a symplectic
submanifold of dimension two does not change the diffeotype, while
blowing up at a point gives us an exceptional divisor of
self-intersection number $-1$. So $b^+=1$ for all symplectic
quotients taken at a regular value.

However, when the action is not semi-free, there is a glitch in the
above argument since in this case the symplectic quotients taken at
regular values are orbifolds,  which causes some technical
difficulties. One might want to use the results established in
\cite{Go00} and compute the change in $b^+$ when passing a critical
value bare-handedly. However, in this paper, we circumvent this by
resorting to the wall crossing formula for the signature of
symplectic quotients developed by Metzler \cite{Me00}, which holds
for Hamiltonian torus actions in general.

This paper is organized as follows. Section \ref{prep} reviews some
basic concepts and results  in symplectic geometry  to set up the
stage. Section \ref{construction} proves for any closed symplectic
four manifold $N$ with $b^+>1$, there exists a Hamiltonian manifold
fibred over $N$ such that the Duistermaat-Heckman function is
non-log-concave.  Section \ref{main examples} applies these results
to construct simply connected examples with the Hard Lefschetz
property. Section \ref{log-concavity} proves the log-concavity
conjecture for Hamiltonian circle actions on six manifolds whose
fixed points sets are either of codimension at most two or only
having four dimensional components with $b^+=1$.

    \,\,\,\,\,\,\,\,\,\,\,\,\,\,\,\,\,

 \textbf{Acknowledgement}: I would like to thank Reyer
Sjamaar for first introducing me to the literature related to the
log-concavity property of the Duistermaat-Heckman measure when I was
a graduate student working with him, and for continuous
encouragement and moral support over the years. I would like to
thank Yael Karshon for a useful remark on an early version of
Theorem \ref{b_2+=1} in this article, and I would like to thank Lisa
Jeffrey for helpful discussions.

%\begin{theorem} \label{Gompf}
%\footnote{The first assertion of Theorem
% \ref{Gompf} is contained in the
  %   statement of Theorem 4.1 of \cite{Gm95}; the second
   %  assertion follows from the discussion following the proof of
    % Observation 7.4 in the same paper.}(\cite{Gm95})
 %    Let $G$ be any
% finitely presentable group. Then there is a closed, symplectic
% $4$-manifold $(M, \omega)$ such that \begin{enumerate}
% \item [(1)] $\pi _1(M)=G$, \item [(2)] The Lefschetz map
% $L_{[\omega]}: H^1(M) \rightarrow H^3(M)$ is trivial.\end{enumerate}
% \end{theorem}

    \,\,\,\,\,\,\,\,\,\,\,\,\,\,\,\,\,

   %\textbf{Acknowledgement}:

%main text

\section{Preliminaries}\label{prep}

\subsection{Intersection form of $4n$ dimensional symplectic
manifolds}

 For a compact orientable manifold $N$ of dimension $4n$,  the intersection form $Q$ on the
$2n$-th integer cohomology of $N$ is a symmetric bilinear form defined by:
\[ Q:H^{2n}(N,\Z) \times H^{2n}(N,\Z) \rightarrow \Z, \,\,\,\, Q([\alpha],[\beta])=<[\alpha] \cup
[\beta] ,[N]>,\] where $[N]$ is the fundamental class of the
manifold $N$.  Using the De Rham model, the corresponding form on
the real cohomology can be defined by:
\[ Q:H^{2n}(N,\R)\times H^{2n}(N,\R)\rightarrow \R, \,\,\,\,
Q([\alpha],[\beta])=\int_N\alpha \wedge \beta.\]

By the Poincar\'e duality, this is a non-degenerate symmetric
bilinear form. We define $b^+$ and $b^-$ to be the dimensions of
maximal positive and negative subspaces of the form, and define the
signature of the manifold $N$ to be $\sigma(N)=b^+-b^-$. Note that
when $N$ is a symplectic manifold, we will assume that the
orientation on $N$ is the one induced by the symplectic form.
%And
%when $N$ is a point, we will adopt the convention that its signature
%equals one.

The following simple looking lemma is actually a key point for our
construction of Hamiltonian manifolds with non-log-concave
Duistermaat-Heckman measure in Theorem \ref{main-construction}.

\begin{lemma}\label{integral-class} Let $(N, \omega_0)$ be a compact symplectic four manifold  such that $b^+>1$ and such that
 $[\omega_0]$ is a rational cohomology
class in $H^2(N)$. And let $Q$ be the intersection form of
$(N,\omega_0)$. Then there exists an integral cohomology class $[c]
\in H^2(N, \R)$ such that
 \begin{equation} \label{tech-condition}
 Q([c],[c])>0 \text{ and }
 Q([c],[\omega_0])=0. \end{equation}
\end{lemma}

\begin{proof} Write $\alpha_1=[\omega_0]$. Since $b^+>1$ and
since $Q(\alpha_1,\alpha_1)>0$, over the field of rational numbers
there exists a basis $\alpha_1,\alpha_2,\cdots, \alpha_r$ of
$H^2(N,\Q)$ such that \begin{itemize} \item [a)] if $ 1\leq i=j\leq
b^+$, then $Q(\alpha_i,\alpha_j)>0$ ; \item [b)]  if $ b^+ < i=j
\leq r$ , then $Q(\alpha_i,\alpha_j)<0$; \item [c)] if $i\neq j$,
then $Q(\alpha_i,\alpha_j) =0$. \end{itemize}

Choose $[c]=\alpha_2 \in H^2(N, \Q)$. Then it is easy to see that
$Q([c],[c])>0$ and $  Q([c],[\omega_0])=0$. Replace $[c]$ by $[nc]$
for some appropriate positive integer if necessary, we get an
integral class $[c]$ such that the condition (\ref{tech-condition})
holds.

\end{proof}
The following non-trivial fact was proved by  Baldridge \cite{Ba04}
using Seiberg-Witten invariants and provides a useful criterion when
a symplectic four manifold must have $b^+=1$.

\begin{theorem} \label{b_2+fourmanifold}(\cite{Ba04}) A symplectic $4$-manifold which admits a circle
action with fixed points must have $b^+=1$.

\end{theorem}

\vskip 5pt

%Now suppose that a $k$\--dimensional torus $T$ acts on $M$ in a
% Hamiltonian fashion and $X$ is the fixed points submanifold of the
% action of $T$. Choose a $T$\--invariant Hermitian inner product on
% $E$ and linearize the action of $T$ on $E_{\delta}$. Then $E$
% decomposes into a direct sum of Hermitian vector bundles $E_{1}
% \oplus E_2\oplus \cdots \oplus E_k$ according to the weights of $T$
% action on $E$.  The structure group of the Hermitian vector bundle
% $E$ reduces to $H:=U(l_1)\times U(l_2)\times \cdots \times U(l_k)
% \subset U(n)$, where $l_i=\text{dim} E_i$, and we can assume that
% $T$ acts on $E_{\delta}$ through some homomorphism $T \rightarrow H$
% such that the image of $T$ lies in the center of $H$. Set $P$ to be
%  the principal $H$ bundle associated to $E=E_1\oplus E_2\oplus \cdots
% \oplus E_k$ and perform the reduction at the zero level with respect
% to the the diagonal Hamiltonian $H$ action on $V^*P\times L$. Then
% the  reduced space is just $E(\epsilon)$.  Since the Hamiltonian $T$
% on $\C^n$  action commutes with the action of $H$, it descends to a
% Hamiltonian action on $E(\epsilon)$. Therefore the  $T$ action
% extends to a Hamiltonian action on the blow up $\widehat{M}$.

\subsection{Duistermaat-Heckman function }

Consider the Hamiltonian action of a torus $T$ on a symplectic
manifold $(M, \omega)$. Let $a\in \frak{t}^*$ be a regular value of
the moment map $\Phi: M \rightarrow \frak{t}^*$. When the action of
$T$ on $M$ is not quasi-free,\footnote{ An action of a group $G$ on a manifold $M$ is called quasi-free if the stabilizers
 of the points are connected.} the quotient $M_a=\Phi^{-1}(a)/T$
taken at a regular value of the moment map is not a smooth manifold
in general. However, the singularity is mild and $M_a$ does admit an
orbifold structure in the sense of Satake \cite{S56}, \cite{S57}. Orbifolds,
although not necessarily smooth, do carry differential structures
such as differential forms, fiber bundles, etc. So the usual
definition of symplectic structures extends to the orbifold case. In
particular, the restriction of the symplectic form $\omega$ to the
level set $\Phi^{-1}(a)$ descends to a symplectic form $\omega_a$ on
the reduced space $M_a$ \cite{A. W. orbifolds}. For our purpose, it
is also important to note that any orbifold is a rational homology
manifold \cite{Ful93} and does satisfy the Poincar\'e duality. Thus
any orbiford has a well defined signature just as in the manifold
case.

The theorems that we are going to state in the rest of this section hold for
general torus actions which are not necessarily semi-free. Their statements
actually involve some basic orbifold related notions, such as a principal bundle over
an orbifold and a diffeomorphism of orbifolds. For basic notions in orbifold theory, we refer
to \cite{R01}, \cite{CR01} and \cite{ALR2007}. For a modern treatment of orbifolds from the viewpoint of Lie groupoids, we refer to
\cite{MO02}. For the foundation of Hamiltonian actions on symplectic orbifolds, the interested reader may consult \cite{LeTo97}.

The following Duistermaat and Heckman theorem \cite{DH82} is a
fundamental result in symplectic geometry.
 \begin{theorem}\label{DH-func}(\cite{DH82}) Consider the  effective Hamiltonian action of a
 $k$ dimensional  torus $T$ on a connected compact $2n$ dimensional symplectic manifold $M$ with moment map $\Phi: M \rightarrow t^*$.
We have that
\begin{itemize} \item [a)]
  at a regular value $a \in
t^*$ of $\Phi$,  the Duistermaat-Heckman function $f$ is computed by
 the following formula:
\[ f(a) =\int_{M_{a}}\dfrac{\omega_{a}^{n-k}}{(n-k)!}\,\,\,\,,\]
where $M_{a}=\Phi^{-1}(a)/T$ is the symplectic quotient,
$\omega_{a}$ is the corresponding reduced symplectic form, and
$M_{a}$ has been given the orientation of $\omega_a^{n-k}$.

\item [b)] if $a, a_0 \in t^*$ lie in the same connected
component $C$ of the regular values of the moment map $\Phi$, then
the reduced space $M_{a}=\Phi^{-1}(a)/T$ is diffeomorphic to
$M_{a_0}=\Phi^{-1}(a_0)/T$; furthermore, let $\Gamma$ be the finite
subgroup of $\, T$ generated by all the finite stabilizer groups
$T_z$, where $z \in \Phi^{-1}(a_0)$, and let
$Z_0=\Phi^{-1}(a_0)/\Gamma$, then using the diffeomorphism  $M_a
\rightarrow M_{a_0}$, the reduced symplectic form on $M_{a}$ can be
identified with
\begin{equation}\label{variation}
\omega_{a}=\omega_{a_0}+<c,a-a_0>,
\end{equation}
where $c \in \Omega^2(M,t^*)$ is a closed $t$-valued two form which
represents the Chern class of the principal torus $T/\Gamma$-bundle
$\pi: Z_0 \rightarrow M_{a_0}$.

\end{itemize}

 \end{theorem}

By the Atiyah-Guillemin-Sternberg convexity theorem (cf. \cite{A82}
and \cite{GS2}), the image of the moment map $\Delta=\Phi(M)$ is a
convex polytope. In fact, $\Delta$ is a union of subpolytopes with
the property that the interiors of the subpolytopes are disjoint and
constitute the set of regular values of $\Phi$. \cite{GLS88} gave an
explicit formula computing the jump in the Duistermaat-Heckman
function $f$ across the wall of $\Phi(M)$. Making use of it, Graham
established the following result, c.f., \cite[Section 3]{Gr96}.

\begin{Proposition} ( \cite{Gr96} ) \label{jump-formula} Suppose $\Phi: M\rightarrow t^*$ is the moment map
of the effective Hamiltonian action of torus $T$ on a connected
compact symplectic manifold $M$. Let $a$ be a point on a codimension
one interior wall of $\Phi(M)$, and let $v\in t^*$ be such that the
line segment $\{a+tv\}$ is transverse to the wall.  For $t$ in a
small open interval near $0$, write $g(t)=f(a+tv)$, where $f$ is the
Duistermaat-Heckman function. Then we have $g_+'(0)\leq g'_-(0)$.
\end{Proposition}

\subsection{ The wall crossing formula for the signature of symplectic
quotients}

Consider the Hamiltonian action of $S^1$ on a compact symplectic manifold $M$ with moment map
$\Phi: M\rightarrow \R$. Let $a< a_1$ be two points in the image of
the moment map such that $a_0$ is the unique critical value between
$a$ and $a_1$. Let $X$ be the set of the critical points of
the moment map $\Phi$ which lies inside $\phi^{-1}(a_0)$. Then each connected component of $X$ is a submanifold of
 $M$, which we will call critical submanifolds of $M$. Let $X_1, X_2,\cdots,X_k$ be all the critical
submanifolds sitting inside $\Phi^{-1}(a_0)$, and let $E_i
\rightarrow X_i$ be the symplectic normal bundle of $X_i$ in $M$.
Then for each $ 1\leq i \leq k$ the Hessian of $\Phi$ gives us a
splitting
\[ E_i =E_i^+ \oplus E_i^- \] of $E_i$ into a direct sum of
positive and negative normal bundles. We denote by $2b_i$  and
$2f_i$ the real dimensions of $E_i^{-}$ and $E_i^{+}$ respectively.
The following Theorem \ref{signature-wall-crossing} of Metzler
\cite{Me00} computes the change in the signature and Poincar\'e
polynomial of symplectic quotients across the critical value $a_0$.
By the way, given a topological space $Y$, throughout this paper we
will always denote by $P(Y)(t)$ its Poincar\'e polynomial.

\begin{theorem} (\cite[pp. 3502, 3518]{Me00})\label{signature-wall-crossing}
Denote the half rank of the symplectic normal bundle $E_i$ by $q_i$.
Then \[
\begin{split} & \sigma(M_{a_1})-\sigma(M_{a})=\sum_{1\leq i\leq
k, q_i \text{ odd}} (-1)^{b_i}\sigma(X_i),\\&
P(M_{a_1})(t)-P(M_{a})(t)=\sum_{i=1}^kP(X_i)(t)\dfrac{t^{2b_i}-t^{2f_i}}{1-t^2},\end{split}
\]

where $M_{a_1}$ and $M_a$ denote the symplectic quotients of the
Hamiltonian $S^1$ action taken at $a_1$ and $a$ respectively.
\end{theorem}

We will also need the following result \cite[Thm. 2.8]{Me00} which
is the orbifold version of a result of Chern, Hirzebruch, and Serre
\cite{CHS57}.

\begin{theorem}\label{fiber-of-signature}Let $P \rightarrow B$ be
a fibre bundle over  $B$ with fibre $F$ such that
\begin{itemize} \item [1)] $P,B,F$ are compact connected oriented
orbifolds; \item [2)] the structure group of $P$ is compact and
connected. \end{itemize}

 If $P,B,F$ are oriented coherently, then $\sigma(P)=\sigma(B) \sigma (F)$.

\end{theorem}

\section{Main construction }  \label{construction}
\label{main-theorem}

A measure defined on $\R^k$ is said to be strictly non-log-concave if it is absolutely
continuous with respect to the Lebesgue measure on $\R^k$ and if the logarithm of its density function \footnote{ Strictly speaking, the density function is well defined only if we declare the  normalization of the Lebesgue measure.} is
a strictly convex function.

\begin{theorem} \label{main-construction} Assume that $N$ is a closed symplectic four
manifold with $b^+>1$, then there exists a
sphere bundle $ \pi: M
 \rightarrow N$ such that there is a symplectic form $\omega$ on $M$
and a Hamiltonian $S^1$ action on $(M,\omega)$  for which the
Duistermaat-Heckman measure is strictly non-log-concave. %\begin{itemize} % \item [a)] Then there exists a sphere
% bundle $ \pi: M
 %\rightarrow N$ such that there is a symplectic form $\omega$ on $M$
%and a Hamiltonian $S^1$ action on $(M,\omega)$  for which the
%Duistermaat-Heckman function is log-concave.

  %\end{itemize}

\end{theorem}

\begin{proof} Note that any
closed symplectic manifold admits another symplectic form whose
cohomology class is integral, see for e.g. \cite[pp.561]{Gom95}.
Without the loss of generality, henceforth we will assume that
$N$ is equipped with a symplectic form $\omega_0$ such that $[\omega_0]$ lies in the
image of $H^2(N,\Z)$ in $H^2(N,\R)$.

 % First we will show for a given integral
% cohomology class $[c]$ in the image of $H^2(N,\Z)$ in $H^2(N,\R)$, there exists a sphere bundle $
% \pi: M \rightarrow N$ such that there is a symplectic form $\omega$
% on $M$ and a Hamiltonian $S^1$ action on $(M,\omega)$.

By Lemma \ref{integral-class}, there exists an integral class $[c]
\in H^2(N, \R)$ which satisfies the condition
(\ref{tech-condition}). Let $\pi_P: P \rightarrow N $ be the principal $S^1$ bundle with
Euler class $[c]$, let $\Theta$ be the connection $1$-form such that
$d\Theta =\pi^*_P c$, let $S^1$ act on $S^2$ by rotation and let $M$
be the associated bundle $P\times_{S^1} S^2$.  Note that the action
of $S^1$ on $S^2$ preserves the standard symplectic form, i.e., the
area form, on $S^2$, and is Hamiltonian. Indeed, in cylindrical
polar coordinates $(\theta, h)$ away from the poles, $0\leq \theta
<2\pi, -1\leq h \leq 1$, the area form $\sigma$ on $S^2$ can be
written as $ d\theta \wedge d h$ and the moment map for the rotating
action of $S^1$ is just the height function $\mu=h$. Since $S^1$ is
the structure group, $\pi: M \rightarrow N$ is a symplectic
fibration over the compact symplectic four-manifold $N$. The
symplectic form $\sigma$ on $S^2$ gives rise to a symplectic form
$\sigma_x$ on each fibre $\pi^{-1}(x)$, $x\in N$; moreover, the
$S^1$-action on $S^2$ induces a fibrewise $S^1$ action on $M$.

Next, we resort to minimal coupling construction to get a closed two
form $\eta$ on $M$ which restricts to the forms $\sigma_x$ on the
fibres. Let us give a sketch of this construction here and refer to
\cite{A. Weinstein} and \cite{GS84} for technical details.

  Consider the closed two form $-d(x \Theta)=
-xd\Theta-dx \wedge \Theta$ defined on $P \times \R$, where $x$ is the linear coordinate on $\R$. It is easy to
see the $S^1$ action on $P \times \R$ given by
$$\lambda (p, x)=(\lambda p, x)$$ is Hamiltonian with moment map
$x$. Thus the diagonal action of $S^1$ on $(P \times \R) \times S^2$
is also Hamiltonian, and $M$ is just the reduced space of $(P\times
\R)\times S^2$ at the zero level. Moreover, the closed two form
$\left(-d(x\Theta) +\sigma \right) \mid _{\text{zero level}}$
descends to a closed two form $\eta$ on $M$ with the desired
property; whereas $x \mid _{\text{ zero level}}$ descends to a globally defined function $H$ on $M$ whose restriction
to each fiber $S^2$ is just the height function $h$.

It is useful to have the following explicit description of $\eta$.
Observe that $d\theta-\Theta$ is a basic form on $(P \times \R)
 \times S^2$. Its restriction to the zero level of $(P \times \R)
 \times S^2$ descends to a one form $\tilde{\theta}$ on $M$ whose
 restriction to each fibre $S^2$ is just $d\theta$. It is easy to see
that on the associated bundle $P\times_{S^1}(S^2-\text{\{two
poles\}})$ we actually have $\eta = H\pi^*_M c+ \tilde{\theta}\wedge dH$.

% For any real number $t_0$, note that the restriction of $\eta -t_0
% \pi_M ^* c$  to fibres are symplectic forms $\sigma_x$. By an
% argument due to Thurston \cite{MS98},  for sufficiently large
% constant $K>0$ the form $K\pi_M^*\omega_0 -t_0 \pi_M^* c+\eta$ is
% symplectic. Equivalently, define $ \omega=\pi_M^*\omega_0 -\epsilon
% t_0 \pi_M^*c+ \epsilon \eta$ for sufficiently small constant
% $\epsilon>0$. Then $\omega$ is a symplectic form on $M$;
% furthermore, the fibrewise $S^1$-action on $(M, \omega)$ is
% Hamiltonian with the moment map $H: M \rightarrow \R$.

Note that the restriction of $\eta$ to each fiber are symplectic
forms $\sigma_x$. Thus by a famous argument due to Thurston
\cite{MS98}, for sufficiently small constant $\epsilon>0$ the form
$\pi_M^*\omega_0 +\epsilon\eta$ is symplectic. Without the loss of generality, we will
assume that $\epsilon$ is so small that
\begin{equation}\label{epsilon} Q([\omega_0],[\omega_0])
 >  \epsilon^2Q([c],[c]).\end{equation}
Having chosen such a symplectic form $\pi^*\omega_0+\epsilon \eta$
on $M$, a simple calculation shows that the fibrewise $S^1$-action
on $(M, \omega)$ is Hamiltonian with the moment map $t:=\epsilon H:
M \rightarrow \R$.

Now let us compute the Duistermaat-Heckman function $f$ . Observe
that at the level set $ -\epsilon < t <\epsilon$, the symplectic
quotient is just $N$ with the reduced symplectic form $\omega_0+t
c$. Therefore by Proposition \ref{DH-func} the Duistermaat-Heckman
function
\[ \begin{split} & f(t)= \int_N \dfrac{1}{2} (\omega_0+t c)^2 \\&=
\dfrac{1}{2}\left( Q([c],[c]) t^2+ 2 Q([c],[\omega_0]) t+
Q([\omega_0],[\omega_0]) \right)\end{split}\]

First note that $\ln f$ is strictly non-log-concave if and only if
$f''f-(f')^2 > 0$. However, a simple calculation shows that
\[ \begin{split} 2\left(f''f-(f')^2
\right)& = \left( Q([c],[c])Q([\omega_0],[\omega_0])
-2Q^2([c],[\omega_0])\right) \\& -\,\,\,\,\,\,   \left(
Q^2([c],[c])t^2+2Q([c],[c])Q([c],[\omega_0])t\right).
\end{split}\]
Since by construction $[c]$ satisfies the condition (\ref{tech-condition}), we have
\[ \begin{split} 2\left(f''f-(f')^2
\right)& =  Q([c],[c])Q([\omega_0],[\omega_0])
 - Q^2([c],[c])t^2.
\end{split}\]

Now it follows easily from  the inequality (\ref{epsilon}) and $ -\epsilon < t < \epsilon$ that
$f''f-(f')^2>0$. Hence the Duistermaat-Heckman function $f$ is
strictly non-log-concave. This finishes the proof of Theorem \ref{main-construction}.

% To prove the assertion (a), let us choose $[c]=[\omega_0]$. Then  we
% have
% \[ \begin{split} 2 \left((f'')^2f-(f')^2 \right)& =
%-Q^2([\omega_0],[\omega_0])-Q^2([\omega_0],[\omega_0])t^2-2Q^2([\omega_0],[\omega_0])
%t \\& =-Q^2([\omega_0],[\omega_0])\left( t+1\right)^2 \leq 0
%\end{split}\]

%This establishes the assertion (a).

\end{proof}

%\begin{example}   Choose $G$ to be any non-K\"ahler group, i.e., a
%finitely presentable group which can not be realized as the
%fundamental group of a compact K\"ahler manifold (see for e.g.,
%\cite{AB96}), and let $N$ be a closed symplectic four manifold with
%$\pi_1(N)=G$. Applying the assertion (a) in Theorem
%\ref{main-construction}, we get a Hamiltonian six manifold $M$ for
%which the Duistermaat-Heckman function is log-concave. Moreover, as
%we explained in the above paragraph, $\pi_1(M)=\pi_1(N)=G$. So $M$
%is not a K\"ahler manifold.

%\end{example}

\begin{remark}  For instance, $K3$ surfaces and complex tori are K\"ahler surfaces with $b^+=3$.
Applying Theorem \ref{main-construction} to the case $N=T^4$,  we recover Karshon's example of a non-log-concave Duistermaat-Heckman measure.
In the general symplectic category, we note that there are many examples of symplectic four manifolds
with $b^+>1$, c.f.,  \cite{Gom95} and \cite{PS00}. Thus Theorem \ref{main-construction} provides us
 with a large family of Hamiltonian manifolds with a non-log concave
 Duistermaat-Heckman measure.

 % However, the
 %assertion (a) there is also interesting on its own. Recall that a well known result of Gompf \cite[Thm. 4.1]{Gom95}
% asserts that the existence of symplectic  structures imposes no
 %restrictions on the fundamental group. Combine this result with
 %the assertion (a), one sees immediately that the existence of a
 %Hamiltonian circle action with log-concave function imposes no
% restrictions on the fundamental group of the underlying manifold.

\end{remark}

\section{Simply connected examples with the Hard Lefschetz property}
\label{main examples}

A compact symplectic manifold $(M, \omega)$ of dimension $2m$ is said to
have the Hard
    Lefschetz property or equivalently to be a Lefschetz manifold if
    and only if for any $0 \leqslant k \leqslant m$, the Lefschetz type map
    \begin{equation}\label{lefschetz property}
    L^k_{[\omega]} : H^{m-k}(M,\R) \rightarrow H^{m+k}(M,\R),\,\,
    [\alpha] \rightarrow [\alpha \wedge \omega^k]
    \end{equation}
    is an isomorphism.

   The follow result allows us to
   construct simply connected Hamiltonian manifolds with the Hard Lefschetz property
   which have a strictly non-log-concave Duistermaat-Heckman measure.

\begin{theorem} \label{Hard-Lefschetz} Assume that $N$ is a closed symplectic four
manifold which satisfies the Hard Lefschetz property and $b^+>1$. Then there exists a
sphere bundle $ \pi: M
 \rightarrow N$ such that there is a symplectic form $\omega$ on $M$
and a Hamiltonian $S^1$ action on $(M,\omega)$  for which the
Duistermaat-Heckman function is strictly non-log-concave; moreover,
$(M, \omega)$ is a compact symplectic manifold which satisfies the Hard Lefschetz property.
\end{theorem}
%\begin{itemize}
%\item [a)] Then there exists a sphere bundle $ \pi: M
 %\rightarrow N$ such that there is a symplectic form $\omega$ on $M$
%and a Hamiltonian $S^1$ action on $(M,\omega)$  for which the
% Duistermaat-Heckman function is log-concave; moreover, we have that
%$(M, \omega)$ satisfies the Hard Lefschetz property.

\begin{proof} Let $(M, \omega)$ be the symplectic
six manifold constructed in the proof of Theorem
\ref{main-construction}. It suffices to show that in the proof of
Theorem \ref{main-construction}, for sufficiently small constant
$\epsilon>0$,  $(M, \pi^*\omega_0 +\epsilon \eta)$ satisfies the
Hard Lefschetz property. The proof  is very similar to the one given
in the proof of  \cite[Prop. 4.2]{Lin07}, though the assumption here
is slightly different from the one used in \cite[Prop. 4.2]{Lin07}.

 Note that the restriction of the cohomology class of $\eta$ to each fibre
$S^2$ generates the second cohomology group $H^2(S^2,\R)$. By the
Leray-Hirsch theorem (see e.g., \cite[Thm.
5.11]{bott-tu;differential-forms}), the pullback map $\pi^*:H^*(N)
\rightarrow H^*(M)$ embeds $H^*(N)$ into $H^*(M)$, and additively
$H^*(M,\R)$ is a free module generated by $1$ and $[\eta]$. It
follows that $[\eta^2]=[\pi^*\beta_2 \wedge \eta]+[\pi^*\beta_4]$
for some closed forms $\beta_2$ and $\beta_4$ on $N$ of degree two
and four respectively.

Choose an $\epsilon>0$ which is sufficiently small such that
 the determinant of the linear
map $L_{[2\omega_0 + \epsilon \beta_2]}: H^1(N,\R) \rightarrow
H^3(N,\R)$ is non-zero and such that \begin{equation}  \label{neq1}
[\omega_{0}]^2 \neq -\epsilon^2 [\beta_4] + \epsilon [\omega_0
\wedge \beta_2].
\end{equation}

 We claim for the $\epsilon$ chosen above, the
symplectic manifold $(M, \pi^*\omega_0 + \epsilon \eta)$ satisfy the
Hard Lefschetz property. By the Poincar$\acute{e}$ duality it
suffices to show the two Lefschetz maps
\begin{equation}\label{injective map1}
L^2_{[\omega ]}: H^1(M,\R) \rightarrow H^5(M,\R)
\end{equation} and
\begin{equation}\label{injective map2}
L_{[\omega ]}: H^2(M,\R) \rightarrow H^4(M,\R)
\end{equation}
are injective. We will give a proof in two steps below.
\begin{enumerate}

\item Since by the Leray-Hirsch theorem $H^1(N)
\xrightarrow [\pi^{\ast}]{\backsimeq} H^1(M)$, to show Map
(\ref{injective map1}) is injective we need only to show that for
any $ [\lambda] \in H^1(N,\R)$ if $L^2_{[\omega ]}
(\pi^{\ast}[\lambda])=0$, then $[\lambda]=0$. A straightforward
calculation shows that
\begin{equation}
\begin{split} 0& = L^2_{[\omega ]} ([\pi^{\ast}\lambda])= \pi^{\ast}\left(2\epsilon
[\omega_0] + \epsilon^2 [\beta_2]\right) \wedge[\pi^{\ast}
\lambda]\wedge [\eta].\end{split}\end{equation} Since  $H^*(M)$ is
free over $1$ and $[\eta]$,  $ 0 = \pi^{\ast}\left(
([2\epsilon\omega_0+ \epsilon^2 \beta_2]) \wedge [\lambda]\right)$.
By our choice of $\epsilon$, the determinant of the linear map
$L_{[2\epsilon\omega_0 + \epsilon^2 \beta_2]}: H^1(N) \rightarrow
H^3(N,\R)$ is non-zero and so we have $[\lambda]=0$.

\item To show that Map
(\ref{injective map2}) is injective it suffices to show that if
$L_{[\omega]}(\pi^{\ast}[\varphi]+ k [\eta])=0$ for arbitrarily
chosen scalar $k$ and second cohomology class $[\varphi] \in
H^2(N,\R)$, then we have $[\varphi]=0$ and $k=0$. Since
$\omega=\pi^{\ast}\omega_0 +\epsilon \eta$ and $[\eta^2]
=[\pi^{\ast}\beta_2 \wedge \eta]+ [\pi^{\ast}\beta_4]$, we have
\begin{equation} \begin{split} 0&=L_{[\omega]}(\pi^{\ast}[\varphi]+ k
[\eta])\\&=\left( \pi^{\ast}[\omega_0\wedge\varphi] + \epsilon k
\pi^{\ast}[\beta_4] \right)+ \left( k \pi^{\ast}[\omega_0 ] +
\epsilon \pi^{\ast}[\varphi] +\epsilon k \pi^{\ast}[\beta_2]
\right)\wedge \eta \end{split}
\end{equation} Since $H(M)$ is a free
module over $1$ and $[\eta]$,  we have that
\begin{equation} \label{part1}\pi^{\ast}[\omega_0\wedge\varphi] + \epsilon k \pi^{\ast}[\beta_4]=0
\end{equation}
\begin{equation}\label{part2}
k \pi^{\ast}[\omega_0 ] + \epsilon \pi^{\ast}[\varphi] +\epsilon k
\pi^{\ast}[\beta_2]=0
\end{equation}
If $k=0$, it follows easily from Equation (\ref{part2}) that
$[\varphi]=0$.  Assume $k \neq 0$, substitute $ \pi^{\ast}[\varphi]=
- \dfrac{1}{\epsilon} k \pi^{\ast}[\omega_0]-k \pi^{\ast}[\beta_2]$
into Equation (\ref{part1}). As a result,
$$ \pi^{\ast}[\omega_0]\wedge (-k \pi^{\ast}[\omega_0]-\epsilon k \pi^{\ast}[\beta_2])+
\epsilon^2 k \pi^{\ast}[\beta_4]=0
$$
For $k\neq 0$, we get
$$\pi^{\ast}([\omega_0])^2 = -\epsilon^2\pi^{\ast}
[\beta_4] + \epsilon \pi^{\ast}[(\omega_0) \wedge \beta_2],$$ which
clearly contradicts Equation (\ref{neq1}).
\end{enumerate}

\end{proof}

\begin{example} Choose any simply connected compact symplectic four
manifold $N$ such that $b^+>1$. (Examples of such  symplectic
manifolds are abundant. For instance, choose $N$ to be $3CP^2 \#
19\overline{C p^2}$, c.f., \cite{Gom95}.) Note that by the Poincar\'e duality any simply connected compact symplectic four manifold
satisfies the Hard Lefscehtz property. Applying Theorem
\ref{Hard-Lefschetz} we get a six dimensional Hamiltonian $S^1$
manifold $(M, \omega)$ which satisfies the Hard Lefschetz property
and which has a strictly non-log-concave Duistermaat-Heckman function. It
then follows easily from the long exact sequence of homotopy groups
for an $S^2$ fibration that $M$ is simply connected as well. The
Hamiltonian manifold $(M, \omega)$ does not admit an $S^1$ invariant
K\"ahler structure since its Duistermaat-Heckman function is
non-log-concave, c.f., \cite{Gr96}.
\end{example}

\section{The log-concavity for torus actions of complexity two}
\label{log-concavity}

The Hamiltonian action of a $k$-dimensional torus on a $2n$
dimensional symplectic manifold is said to be of complexity two if
$n-k=2$. Theorem \ref{b_2+=1} gives a useful criterion to ensure the
log concavity of the Duistermaat-Heckman measure for a Hamiltonian
torus action of complexity two.

\begin{theorem}\label{b_2+=1} Assume that the action of  a torus  $T$  on
a connected compact symplectic manifold $M$ is an effective
Hamiltonian action of complexity two with  moment map $\Phi:
M\rightarrow t^*$. And assume that for any regular value $\xi \in
t^*$ of $\Phi$, the symplectic reduced space
$M_{\xi}=\Phi^{-1}(\xi)/T$ has that $b^+=1$. Then the
Duistermaat-Heckman function $f$ is log-concave.
\end{theorem}

\begin{proof} In view of Proposition \ref{jump-formula}, to establish the log concavity
of $f$ on $\Phi(M)$, it suffices to show that the restriction of $
\ln f$ to each connected component of the set of regular values of
$\Phi$ is concave. Let $C$ be such a component, let $v \in
\frak{t}^*$, and let $\{a+tv\}$ be a line segment in $C$ passing
through a point $a \in C$, where the parameter $t$ lies in some
small interval containing $0$. We need to show that $g(t):=f(a+tv)$
is log-concave, or equivalently, $g''g-(g')^2\leq 0$.

  It follows from Theorem
\ref{DH-func} that the Duistermaat-Heckman function $f$ is computed by
\begin{equation} \begin{split} f(a+tv) &=\dfrac{1}{2}
\int_{M_a}\left(
\omega_a+tc\right)^2\\&=\dfrac{1}{2}\left(Q([c],[c])t^2+2Q([c],[\omega_a])t+Q([\omega_a],[\omega_a])\right),
 \end{split}
\end{equation}where $M_a=\Phi^{-1}(a)/T$ is the reduced space at $a \in C$, $\omega_a$ is the reduced symplectic form on it,
and $c \in \Omega^2(M_a)$ is a closed two form depending only on $v$
in $C$. Consequently,
\[ \begin{split} 2(g''g-(g')^2)&=\left( Q([c],[c])Q([\omega_a],[\omega_a])
 -2Q^2([c],[\omega_a])\right) \\&\,\,\,\,\,- \left( Q^2([c],[c])
t^2+2Q([c],[c])Q([c],[\omega_a])t\right) \end{split}\]
 Since  $b^+(M_a)=1$, there exists a real basis
$\alpha_1, \alpha_2,\cdots,\alpha_k$ of $H^2(N,\R)$ such that
$[\omega_a]=r\alpha_1$ for some positive constant $r$ and such that
\[  Q(\alpha_i,\alpha_j)=    \begin{cases}

1, &\text{ if $i=j=1$ ;} \\
-1, &\text{ if $2 \leq i=j \leq k $;}\\
0,  &\text{ otherwise.} \end{cases} \] Write
$[c]=\sum_{i=1}^k\lambda_i\alpha_i$ for some real scalars
$\lambda_i$. Then we have
\[ \begin{split}\label{qudratic-poly} 2(g''g-(g')^2)&
=-\left((\lambda_1^2-\lambda_2^2-\cdots-\lambda_k^2)^2t^2+2(\lambda_1^2-\lambda_2^2-\cdots-\lambda_k^2)\lambda_1rt\right)
\\&\,\,\,\,\,\,+\left((\lambda_1^2-\lambda_2^2-\cdots-\lambda_k^2)r^2 -2\lambda_1^2r^2\right)\\&=
-(\lambda_1^2-\lambda_2^2-\cdots-\lambda_k^2)^2t^2-2(\lambda_1^2-\lambda_2^2-\cdots-\lambda_k^2)\lambda_1rt
 \\ &\,\,\,\,\,\,-\left(\lambda_1^2+\lambda_2^2+\cdots+\lambda_k^2\right)r^2 .
\end{split}\]
If the leading coefficient
$-(\lambda_1^2-\lambda_2^2-\cdots-\lambda_k^2)^2$ of the above
polynomial equals zero, then obviously we have $g''g'-(g')^2\leq 0$.
Otherwise, $2(g''g'-(g')^2)$ is a quadratic polynomial with a
negative leading coefficient. Furthermore, the discriminant of this
quadratic polynomial is
\[\begin{split}
\Delta&=4(\lambda_1^2-\lambda_2^2-\cdots-\lambda_k^2)^2\lambda_1^2r^2-
4(\lambda_1^2-\lambda_2^2-\cdots-\lambda_k^2)^2(\lambda_1^2+\lambda_2^2+\cdots+\lambda_k^2)r^2
\\&=4(\lambda_1^2-\lambda_2^2-\cdots-\lambda_k^2)^2\left(\lambda_1^2r^2-(\lambda_1^2+\lambda_2^2+\cdots+\lambda_k^2)r^2\right)
\\&=-4(\lambda_1^2-\lambda_2^2-\cdots-\lambda_k^2)^2(\lambda_2^2+\cdots+\lambda_k^2)r^2
\end{split}\]
  which is clearly non-positive. Thus $2(g''g-(g')^2)$ has to be
  negative for all $t \in I$. This finishes the proof of the
  theorem.

\end{proof}

Theorem \ref{log-concavity-conjecture} and Theorem
\ref{log-concavity-conjecture-2} below give a rather satisfactory
answer to the question when the log-concavity conjecture holds for
$S^1$ actions on six manifolds. To prove them, we need to establish
the following key lemma.

\begin{lemma}\label{wall-acrossing-b+} Suppose the action of $S^1$
on a connected compact symplectic six manifold $M$ is Hamiltonian.
 Let $\Phi: M\rightarrow \R$ be the moment map, and
let $M_a:=\Phi^{-1}(a)/S^1$ be the symplectic quotient taken at the
regular value $a$ of $\Phi$.  Then $b^+(M_a)$ remains constant as
$a$ runs through all the regular values of $\Phi$.

\end{lemma}

\begin{proof} Let
$a_0=\text{ minimum }<a_1<\cdots<a_k=\text{ maximum }$ be all the
critical values of $\Phi$.

By Theorem \ref{DH-func} the diffeotype of $M_a$ remains unchanged
on each open interval $(a_{i-1},a_i)$, $ 1\leq i \leq k$. Thus
$b^+(M_a)$ is a constant on each open interval $(a_{i-1},a_i)$. Next
we note that for dimension reasons, if a critical submanifold $X$ is
neither the minimum nor the maximum submanifold, the signature of
the Hessian of $\Phi$ at $X$ can only be of the form $(2,2p)$ or
$(2,2q)$ for some integers $p,q>0$.  It follows that any such
critical submanifold can be of dimension at most $2$. By the way,
when the signature of the Hessian of $\Phi$ at a critical
submanifold $X$ is of the form $(2p,2q)$, we will say that the
critical submanifold $X$ is of type $(2p,2q)$ \footnote{The standard
terminology is to say that the critical submanifold is of signature
$(2p,2q)$. Because in our paper the word "signature" has been
reserved to refer to something else, we use the word "type" here to
avoid the confusion.}.

Now let $X_i$ be all the critical submanifolds sitting inside
$\Phi^{-1}(a_i)$, and let $(2f_i,2b_i)$ be the signature of the
Hessian of $\Phi$ at $X_i$. Then by Theorem
\ref{signature-wall-crossing}, the change in the signature of
symplectic quotients when passing the critical value $a_i$ is
computed by
 \begin{equation} \label{signature-change}\sum_{1\leq i \leq k, q_i\text{ odd}} (-1)^{b_{i}}\sigma(X_i),\end{equation}
while the change in the Poincar\'e polynomial is computed by
\begin{equation} \label{Poincar-polynomial-change}
\sum_{i=1}^kP(X_i)(t)\dfrac{t^{2b_i}-t^{2f_i}}{1-t^2}.\end{equation}

Note that if  the dimension of $X_i$ is two, then the signature of
the Hessian of $\Phi$ at $X_i$ is $(2,2)$ and $\sigma(X_i)=0$. So
$X_i$ does not have any contribution in either Equation
(\ref{signature-change}) or Equation
(\ref{Poincar-polynomial-change}).

Let $N_1$ be the number of the type $(2,4)$ isolated fixed points
sitting inside $\Phi^{-1}(a_i)$, and let $N_2$ be the number of the
type $(4,2)$ isolated fixed points sitting inside $\Phi^{-1}(a_i)$.
Then when passing the critical value $a_i$, the change in the
signature of symplectic quotients equals $N_1-N_2$, whereas the
change in the second Betti number equals $N_2-N_1$. Therefore the
change in the sum $\sigma+b_2$ is null. Note that for any four
manifold $b^+= \dfrac{1}{2}(\sigma+b_2)$. So the change in
$b^+(M_a)$ when $a$ passes through the critical level $a_i$ is also
null. This finishes the proof of the lemma.

\end{proof}

\begin{theorem}\label{log-concavity-conjecture} Let $(M, \omega)$ be a compact connected symplectic six manifold
equipped with a Hamiltonian $S^1$ action whose fixed points set has
codimension greater than or equal to four.  Then the
Duistermaat-Heckman function of $M$ is log-concave.
\end{theorem}

\begin{proof} Let  $\Phi$ be the moment map of the $S^1$ action on
$M$ such that $a_0 \in \R$ is the minimum value, let $F$ be the
unique local minimum fixed points submanifold in $\Phi^{-1}(a_0)$
which is of codimension $k$, and let $E$ be the symplectic normal
bundle of $F$ in $M$. Choose an $S^1$ invariant Hermitian inner
product on $E$ such that $E$ becomes a Hermitian vector bundle.
Denote by $P$ the principal $U(k)$-bundle, i.e., the unitary frame
bundle, associated to $E$ and choose a connection on it. This gives
a projection map $TP \rightarrow VP$, where $TP$ is the tangent
bundle of $P$ and $VP$ is the bundle of vertical tangent vectors.
Dually, we have an embedding $i: V^*P \rightarrow T^*P$. Let
$\omega_P$ be the standard symplectic form on the cotangent bundle
$T^*P$. Then the $U(n)$ action on $P$ lifts to an action on $V^*P$
which is Hamiltonian with respect to the two form $i^*\omega_P$ on
$V^*P$.

Consider the diagonal Hamiltonian action of $U(k)$ on $V^*F \times
 C^k$ and perform reduction at the zero level. Then we get a
closed two form  which is non-degenerate on a tubular neighborhood
$E_{\delta}$ of $F$.   Since the standard $S^1$ action on $C^n$
commutes with the $U(k)$-action, it descends to a Hamiltonian action
on $E_{\delta}$. By the equivariant Darboux theorem, we can identify
the above Hamiltonian $S^1$ manifold $E_{\delta}$ with an $S^1$
invariant open neighborhood of $F$ in $M$.  Then by a reduction by
stage argument, it is easy to see that the  for any $a>a_0$
sufficiently close to the minimum value $a_0$, as a topological
space the symplectic reduced space $M_{a}=\Phi^{-1}(a)/S^1$ is just
a weighted $CP^k$-bundle over $F$. Indeed, when $k=6$,
topologically $M_a$ is a weighted projective space $CP^2$, and
when $k=4$, is a weighted $CP^1$-bundle over the surface $F$. We
claim that in both cases $b^+=1$. In the former case, since the
rational cohomology of the weighted projective space $CP^2$  is
isomorphic to that of the ordinary projective space $CP^2$( c.f.,
\cite[pp. 3500]{Me00}) , we have that $b^{+}=1$. In the latter case,
the restriction of the reduced symplectic form $\omega_a$ on $M_a$
to each fiber, a weighted projective space $CP^1$, generates its
second cohmology which is one dimensional. So it follows easily from
the Leray-Hirsch theorem \footnote{The spectral sequence argument
given in \cite[pp.170]{bott-tu;differential-forms} can be easily
adapted to show that the Leray-Hirsch theorem does extend to this
case.} that $H^2(M_a,\R)$ is two dimensional. Beside, it is easy to
see from Theorem \ref{fiber-of-signature} that the signature of
$M_a$ is zero. Therefore we have $b^+(M_a)=1$. Applying Lemma
\ref{wall-acrossing-b+} we have that all the symplectic quotients
taken at regular values of $\Phi$ satisfy $b^+=1$. By Theorem
\ref{b_2+=1} the Duistermaat-Heckman function of $M$ has to be
log-concave. This finishes the proof of the theorem.

\end{proof}
\begin{remark}  Alternatively, when symplectic quotients are regular manifolds, one can show that $b^+(M_a)=1$ for
$a_0<a<a_1$ using Theorem \ref{b_2+fourmanifold}. It is interesting
to notice that in the case codimension is greater than or equal to
$4$, $M_a$ admits an effective fibrewise $S^1$ action which has
fixed points.\end{remark}

\begin{theorem}\label{log-concavity-conjecture-2} Let $M$ be a compact connected symplectic six manifold
equipped with a Hamiltonian $S^1$ action whose fixed points set has
components of dimension four. Then we have that
\begin{itemize} \item [a)] there are only
two such components of dimension four: the unique minimum
submanifold and the unique maximum submanifold;
\item [b)] $b^+(\text{minimum})=b^+(\text{maximum})$. \end{itemize}
If in addition, we assume
$b^+(\text{minimum})=b^+(\text{maximum})=1$, then the
Duistermaat-Heckman function of the Hamiltonian manifold $M$ is
log-concave.
\end{theorem}

\begin{proof} It is a well known result that the any level set of the moment map is connected, c.f., \cite{A82}. In particular,
$M$ has a unique local minimum and a unique local maximum. Thus if a critical submanifold $F$ is neither minimum
nor maximum, then it must be of signature $(2p,2q)$ for some
integers $p,q>0$. Now Assertion (a) follows easily from this observation. Next using the equivariant Darboux theorem, it is
easy to see that for a regular value $a$ sufficiently close to the
minimum value of the moment map $\Phi$, as a topological space the
symplectic quotient $M_a=\Phi^{-1}(a)/S^1$ can be identified with
the minimum submanifold.  Then it follows from Lemma
\ref{wall-acrossing-b+} that
$b^+(M_a)=b^+(\text{minimum})=b^+(\text{maximum})$ for any regular
value $a$ of $\Phi$. This proves Assertion (b). The last assertion
in the theorem now follows easily from Theorem \ref{b_2+=1}.

\end{proof}

% \appendix

% \section{}
% \label{}

% Bibliographic references with the natbib package:
% Parenthetical: \citep{Bai92} produces (Bailyn 1992).
% Textual: \citet{Bai95} produces Bailyn et al. (1995).
% An affix and part of a reference:
%   \citep[e.g.][Ch. 2]{Bar76}
%   produces (e.g. Barnes et al. 1976, Ch. 2).

\end{document}